\documentclass[12pt]{article}
\usepackage{authblk}
\usepackage[a4paper,margin=1in]{geometry}
\usepackage{graphicx}
\usepackage{amsmath,amssymb}
\usepackage{caption}
\usepackage{booktabs}
\usepackage{hyperref}
\usepackage{cite}
\usepackage{algorithm}
\usepackage{algpseudocode}
\usepackage{subcaption} 
\usepackage{float}
\usepackage{booktabs}
\usepackage{siunitx}
\usepackage[numbers]{natbib}

\title{Simulation Platform To Evaluate Inversion Techniques For Magnetic Resonance Elastography Data}

\author[1*]{Yashasvi Verma}
\author[2]{Jakob Schattenfroh}
\author[2]{Ingolf Sack}
\author[3]{Silvia Budday}
\author[1]{Paul Steinmann}
\author[4]{Luca Heltai}
\affil[1]{Institute of Applied Mechanics, Friedrich-Alexander-Universit\"at, Erlangen-N\"urnberg, Germany}
\affil[2]{Elastography Group, Charit\'e - Universit\"atsmedizin Berlin, Germany}
\affil[3]{Institute of Continuum Mechanics and Biomechanics, Friedrich-Alexander-Universit\"at, Erlangen-N\"urnberg, Germany}
\affil[4]{Numerical Analysis Group, Universit\`{a} di Pisa, Pisa, Italy}
\affil[*]{\textit {yashasvi.verma@fau.de}}

\date{}

\begin{document}

\maketitle

\begin{abstract}
Magnetic Resonance Elastography (MRE) has become an essential tool in assessing the mechanical properties of soft tissues \textit{in-vivo}, prompting significant progress in new sequencing techniques and inversion algorithms. This creates a need for a benchmarking framework to promote uniformity and accessibility. To address this, we introduce a comprehensive \textit{in-silico} dataset acquired by solving the forward Finite Element calculations of shear wave propagation in a linear visco-elastic material. This dataset aims to serve as a platform for evaluating inversion schemes by providing data that can be used as numerical input with known mechanical properties to these methods. It includes simulations on homogeneous cuboidal domains of varying spatial and temporal resolution, and an extension to more physiological variations, including material inhomogeneity and internal arterial pulsation. 
To exhibit the importance of this benchmarking, we present a comprehensive case study using raw simulated data as an input to a direct inversion (DI) scheme, which allows for an expedient local/point-wise inversion of the measured data of wave propagation into the underlying material parameters. When aiming to reconstruct the parameters describing the linear visco-elastic material behavior via DI, we find that due to compromised convergence properties of frequency-domain stencils, stemming from truncation and subtractive cancellation errors, the reconstruction accuracy depends non-monotonically on the spatial and temporal resolution of the measurement grid. For inhomogeneous domains, the reconstruction was successful with notable interface boundaries when the material properties change perpendicular to the wave direction. In the presence of pressurized vascular inclusions, a general stiffening of the domain was noted, as the recovered shear modulus was higher than the one assumed in forward modeling.
Our study highlights the potential of this dataset as a vital benchmarking tool for advancing the development and refinement of MRE techniques, ultimately contributing to more accurate and reliable assessment of soft tissue mechanics.
\end{abstract}


\section{Introduction}\label{introduction}
Magnetic Resonance Elastography (MRE) is a non-invasive technique for assessing the mechanical properties of soft tissues \textit{in-vivo}. These mechanical characteristics can be linked to various clinical conditions; for example, changes in MRE-derived stiffness can indicate the presence or progression of tumors, tissue degeneration, or atrophy \cite{pepin_magnetic_2015, sack_structure-sensitive_2013}. MRE has been widely applied to study the mechanical behavior of organs, such as the liver, kidney, breast, and brain \cite{hirsch_magnetic_2017}. Given its broad range of applications, MRE protocols are often tailored to meet specific diagnostic or research needs. A standard MRE workflow typically consists of: (i) inducing time-harmonic mechanical waves in the tissue, (ii) imaging the resulting cyclic displacements to capture wave propagation, and (iii) reconstructing spatial maps of tissue stiffness within the region of interest.

A variety of quantitative inversion techniques have been developed to reconstruct stiffness maps, or elastograms, from displacement fields \cite{meyer_comparison_2022}. Prior to inversion, the measured wave data typically undergoes pre-processing steps such as phase unwrapping, noise filtering, and curl filtering. Phase unwrapping algorithms are essential to correct for phase discontinuities caused by high-amplitude displacements. Noise is commonly reduced using low-pass filters, while curl-based or high-pass filtering is applied to isolate shear waves from low-frequency longitudinal components. The extraction of mechanical parameters from the processed wave data is based on the mathematical description of shear wave propagation in the medium. Several protocols exist for this purpose, including direct inversion (DI), using measured displacements or iterative parameter identification via forward modeling \cite{honarvar_comparison_2016}. The selection of a particular approach depends on factors such as application requirements, computational resources, and the desired level of accuracy.

The Finite Element Method (FEM) has been a vital tool for understanding various aspects of mechanics and has been adopted in the field of biomechanics. FEM is an instrument to effectively model tissue mechanics, thus supporting a better understanding of the human body. In the present context, FEM can be used to mimic the tissue under MRE loading conditions to further analyze, verify and validate the experimental schemes. Thus, we will here exploit the robustness and adaptability of FEM to provide data for benchmarking inversion schemes.\\ 
Honarvar et. al.\cite{honarvar_comparison_2016} used FEM modeling based on phantom material experiments to draw a comparison between various direct and iterative inversion techniques. \textit{In-silico} studies are a good indicator of visualizing and testing the effects of frequency, voxel size and noise in excitation of homogeneous and heterogeneous materials. The modeling domain can be extended to include complex geometries that are more anatomically accurate  \cite{mcgrath_silico_2021}. The findings from FEM models can be interpolated and extended to \textit{ex-vivo} and \textit{in-vivo} data \cite{honarvar_comparison_2017}. The effects of different parametric aspects can be reported as Hollis et al. \cite{hollis_finite_2016} elaborate on the importance of the size of embedded inclusions on the ability to re-extract mechanical parameters using DI schemes. This indicates the use of MRE to accurately predict the shear modulus of heterogeneous soft tissues with degrees of variability as established by Barnhill et al. \cite{barnhill_heterogeneous_2018} for the brain tissue. The anatomical vagaries and geometries can play a crucial role in shaping the wave pattern and hence the errors in mechanical characterization \cite{mcgrath_magnetic_2016}. The effects of the correlation between wave excitation direction and anisotropic fiber orientation on the elastograms has been shown by Mcgarry et al. \cite{mcgarry_mapping_2022}. The inversion methods can be enhanced using FEM models as basis for testing before experimental validation as shown by Kwon et al.\cite{kwon_shear_2009}. They propose an algorithm to reduce the artifacts in shear modulus maps as generated by noisy data. These results can be symbiotic and studies have been done to improve the Finite Element representation of constitutive models that appropriately depict the soft tissue, and the effects of their parameters on the inversion accuracy \cite{mcgarry_suitability_2015}.

The objective of this study is to provide simulated data on reference domains that can be utilized as a shared resource. These controlled datasets mimic tissue-like structures under MRE loading conditions and can be used as a benchmark standard against which the plethora of inversion schemes can be tested such as those available at \url{https://bioqic-apps.charite.de}. The simulations include basic homogeneous domains as well as those with varying material composition and pulsating blood vessels. This presents a variety of complexity against which post-processing protocols from shear wave isolation to inversion techniques can be uniformly assessed for robustness to variations in experimental conditions. The paper begins with a section on the description of the forward problem. This illustrates the FEM model along with the boundary conditions applied to initiate external excitation and reduce boundary reflection. This section also elucidates the 3D-1D immersed boundary method used to embed pulsating inclusions into the domain to mimic blood vessels deforming under arterial pressure. The following section discusses a DI scheme used for generating elastograms. The displacement field from the simulated data is inputted into this inversion scheme to retrieve material parameters. This procedure serves as a case study on how to use the benchmark dataset. The study also reports in detail on the accuracy of the DI highlighted. 

\section{Methods}
\subsection{In-silico Framework for Data Acquisition}\label{FEmodel}
We use a standard finite element formulation for the elastodynamic wave equation, applied here to a linear viscoelastic solid. The computational domain is a 3D cube (representing tissue) with inhomogeneous, linear, and isotropic material behavior. The governing equations for small deformations are given by the linear momentum balance and by a properly chosen constitutive law. In the absence of body forces, the strong form of the initial-boundary value problem reads: \begin{equation}\label{cauchy_eq}
\begin{aligned}
\rho\ddot{\mathbf{u}}(\mathbf{x},t) - \nabla \cdot \boldsymbol{\sigma}(\mathbf{x},t) &= \bold{0}, && \text{for } \mathbf{x}\in \Omega, \quad t\in[0,T], \\
\mathbf{u}(\mathbf{x},t) &= \mathbf{u}_D(\mathbf{x},t), && \text{for } \mathbf{x}\in \Gamma_D, \\
\boldsymbol{\sigma}(\mathbf{x},t)\cdot \mathbf{n} &= \mathbf{t}_N(\mathbf{x},t), && \text{for } \mathbf{x}\in \Gamma_N, \\
\mathbf{u}(\mathbf{x},0) &= \mathbf{0}, && \text{for } \mathbf{x}\in \Omega\\
\dot{\mathbf{u}}(\mathbf{x},0) &= \mathbf{0}, && \text{for } \mathbf{x}\in \Omega,
\end{aligned}
\end{equation} 
where $\rho$ is the density, $\mathbf{u}(\mathbf{x},t)$ is the displacement field, $\ddot{\mathbf{u}}$ denotes material acceleration, and $\mathbf{n}$ is the outward unit normal on the boundary. The $\mathbf{x}$ and $t$ are dropped from further notation for clarity. $\Gamma_D$ and $\Gamma_N$ are the portions of the boundary with prescribed displacements $\boldsymbol{u}_D$ (Dirichlet conditions) and tractions $\mathbf{t}_N$ (Neumann conditions), respectively. In our simulations we apply time-harmonic Dirichlet excitations (to generate shear waves) and take $\mathbf{t}_N=\mathbf{0}$, since external tractions are absent.

\subsubsection{Constitutive Model}
We adopt a Kelvin–Voigt linear viscoelastic material model, which can be thought of as consisting of an elastic spring (shear modulus $\mu$) and a viscous dashpot (viscosity $\eta$) in parallel. The Cauchy stress $\boldsymbol{\sigma}$ is related to the infinitesimal strain $\boldsymbol{\epsilon}(\mathbf{u})$ and strain-rate $\dot{\boldsymbol{\epsilon}}(\mathbf{u})$ by 
\begin{equation}
    \label{constitutive}
    \boldsymbol{\sigma}(\mathbf{u}) = \mathbb{E} : \boldsymbol{\epsilon}(\mathbf{u}) + \mathbb{N} : \dot{\boldsymbol{\epsilon}}(\mathbf{u}),
\end{equation} 
where $\mathbb{E}$ and $\mathbb{N}$ are fourth-order elasticity and viscosity tensors, respectively. For a linear isotropic solid, $\mathbb{E}$ has the form $\mathbb{E} = \lambda \mathbf{I}\otimes\mathbf{I} + 2\mu \mathbb{I}_s$, with $\lambda$ and $\mu$ being the Lamé constants (elastic bulk and shear moduli) and $\mathbb{I}_s$ the symmetric identity tensor (so that $\mathbb{I}_s:\boldsymbol{\epsilon} = \boldsymbol{\epsilon}$). Similarly, the viscosity tensor can be written as $\mathbb{N} = \lambda_v\mathbf{I}\otimes\mathbf{I} + 2\eta\mathbb{I}_s$, where $\eta$ is the shear viscosity (and $\lambda_v$ a bulk viscosity, here taken zero for simplicity). The small strain tensor is given by $\boldsymbol{\epsilon}(\mathbf{u}) = \frac{1}{2}[\nabla \mathbf{u} + \nabla \mathbf{u}^T]$, or, in index notation, $\epsilon_{ij}(u) = \frac{1}{2}[u_{i,j} + u_{j,i}]$ (with $u_{i,j}=\partial u_i/\partial x_j$), and $\dot{\epsilon}_{ij}$ is its material time derivative. Under the assumption of nearly incompressible material behavior, the first Lam\'{e} parameter $\lambda$ is large and the volumetric stress contribution $\lambda \nabla \cdot \mathbf{u}$ gives rise to a hydrostatic pressure field acting as a Lagrange multiplier enforcing incompressibility. While $\nabla \cdot \mathbf{u}$ is small, the product $\lambda \nabla \cdot \mathbf{u}$ does not vanish in general. In Magnetic Resonance Elastography, longitudinal wave components and spatial variations of the hydrostatic pressure are removed by curl filtering of the measured displacement field \cite{sack_structure-sensitive_2013}. As a result, the remaining displacement field is divergence-free in the filtered sense, and the volumetric stress contributions do not enter the governing equations for the shear-wave motion. Under these assumptions, the constitutive law effectively reduces to
\begin{equation}
    \label{sigma_eq}
    \boldsymbol{\sigma}(\mathbf{u}) = 2\mu \boldsymbol{\epsilon}(\mathbf{u}) + 2\eta \dot{\boldsymbol{\epsilon}}(\mathbf{u}),
\end{equation}
which is the Kelvin-Voigt stress-strain relation with shear modulus $\mu$ and shear viscosity $\eta$. This constitutive model is commonly used to represent soft tissue viscoelasticity in elastography studies.

\subsubsection{Weak form}\label{sub_weakform}
Following the standard Galerkin procedure, we derive the weak (variational) form of the elastodynamic problem. Let $H^1(\Omega)^{\text{dim}}$ be the Sobolev space of functions with square-integrable first derivatives, and let $\Gamma_D$ denote the Dirichlet boundary where prescribed displacements are applied. We define the following spaces:
\begin{equation}
    \label{spaces}
    \begin{aligned}
        V &= \{\mathbf{w} \in H^1(\Omega)^{\text{dim}}\} \\
        V_D &= \{\mathbf{u} \in H^1(\Omega)^{\text{dim}} \mid \mathbf{u} = \mathbf{u}_D \text{ on } \Gamma_D\} \\
        V_0 &= \{\mathbf{v} \in H^1(\Omega)^{\text{dim}} \mid \mathbf{v} = \mathbf{0} \text{ on } \Gamma_D\}
    \end{aligned}
\end{equation} 
where $V_0$ is the space of test functions (virtual displacements) that vanish
on the Dirichlet boundary, while $V_D$ is the affine space of displacements that
match the Dirichlet boundary conditions. We multiply the strong form 
\eqref{cauchy_eq} by an arbitrary $\mathbf{v}\in V_0$, integrate over $\Omega$,
and apply integration by parts (the divergence theorem) on the stress term. This
yields the principle of virtual work (balance of internal, inertial, and
external work): 
\begin{equation}
    \label{weakform_eq}
    \int_{\Omega} \rho\ddot{\mathbf{u}}\cdot \mathbf{v}\,d\Omega + \int_{\Omega} \boldsymbol{\sigma}(\mathbf{u}): \boldsymbol{\epsilon}(\mathbf{v})\,d\Omega = \int_{\Gamma_N} \mathbf{t}_N \cdot \mathbf{v}\,d\Gamma~, \qquad \forall\mathbf{v}\in V_0,
\end{equation} where $\boldsymbol{\epsilon}(\mathbf{v} )= \frac{1}{2}[\nabla \mathbf{v} + \nabla \mathbf{v}^T]$ is the virtual strain. In deriving \eqref{weakform_eq}, we have used $\mathbf{v}=\bold{0}$ on $\Gamma_D$ and the symmetry of $\boldsymbol{\sigma}$. For our case with no applied traction on $\Gamma_N$ ($\mathbf{t}_N=\mathbf{0}$), the right-hand side vanishes. Substituting the constitutive relation \eqref{sigma_eq} into \eqref{weakform_eq}, the weak form becomes: \begin{equation}
    \label{weakform_final}
\int_{\Omega}\rho\ddot{\mathbf{u}}\cdot \mathbf{v}\,d\Omega + \int_{\Omega} 2\mu \boldsymbol{\epsilon}(\mathbf{u}): \boldsymbol{\epsilon}(\mathbf{v})\,d\Omega + \int_{\Omega} 2\eta \dot{\boldsymbol{\epsilon}}(\mathbf{u}): \boldsymbol{\epsilon}(\mathbf{v})\,d\Omega = 0~, \qquad \forall\mathbf{v}\in V_0.
\end{equation} 

Equation \eqref{weakform_final} is the variational form of the elastodynamic
problem where Dirichlet boundary conditions are enforced strongly, which serves
as the starting point for the finite element discretization. It embodies the
balance of linear momentum in a virtual-work sense (D’Alembert principle), with
the first term representing inertial forces and the latter terms the internal
elastic and viscous forces.

\subsubsection{Finite Element Discretization}
\label{sub_FEform}

We discretize the weak form \eqref{weakform_final} in space using FEM. The domain $\Omega$ is partitioned into $N_e$ subdomains
(the finite elements, typically chosen to be hexahedral or tetrahedral
subddomains), and we choose piecewise-polynomial basis (shape) functions
$\mathbf{N}_a(\mathbf{x})$ (with $a$ indexing global nodes and vector
components) in $\Re^{dim}$, that are continuous and have compact support to
generate a finite dimensional subspace $V^h$ of $V$. The trial and test fields
are then approximated as $\mathbf{u}^h(\mathbf{x},t) = \sum_{a}
\mathbf{N}_a(\mathbf{x}) d_a(t)$ and $\mathbf{v}^h(\mathbf{x}) = \sum_{b}
\mathbf{N}_b(\mathbf{x}) w_b$, where $d_a(t)$ are the time-dependent nodal
displacement coefficients and $w_b$ are arbitrary constant test nodal
coefficients. Inserting these expansions into \eqref{weakform_final} and
equating the coefficients of the independent $w_b$ yields the semi-discrete
equation of motion in matrix form:
\begin{equation}\label{FEM_matrix} 
    \mathbf{M}\ddot{\mathbf{d}}(t) + \mathbf{C}\dot{\mathbf{d}}(t) + \mathbf{K}\mathbf{d}(t) = \mathbf{f}(t),
\end{equation}
where $\mathbf{d}(t)$ is the global vector of unknown nodal coefficients. Here,
$\mathbf{M}$, $\mathbf{C}$, and $\mathbf{K}$ denote the global mass, damping (viscous), and stiffness
matrices, respectively, assembled in the standard way from elemental
contributions. The entries of these matrices are given by summing the usual integrals
over an element $e$ (with domain $\Omega_e$): 
\begin{align}
    M_{ab} = \sum_e M_{ab}^{(e)} & := \int_{\Omega} \rho\, \mathbf{N}_a \cdot \mathbf{N}_b \, d\Omega, \\
    K_{ab} = \sum_e K_{ab}^{(e)} & := \int_{\Omega} 2\mu\, \boldsymbol{\epsilon}(\mathbf{N}_a) : \boldsymbol{\epsilon}(\mathbf{N}_b) \, d\Omega, \\
    C_{ab} = \sum_e C_{ab}^{(e)} & := \int_{\Omega} 2\eta\, \boldsymbol{\epsilon}(\mathbf{N}_a) : \boldsymbol{\epsilon}(\mathbf{N}_b) \, d\Omega,
\end{align}
and $\mathbf{f}(t)$ is the global force vector arising from any applied body
forces or Neumann boundary conditions (zero in our case). By assembling all
element matrices, and summing over the local elements, we obtain the global
system \eqref{FEM_matrix}. Equation \eqref{FEM_matrix} is a system of
second-order ordinary differential equations in time, which we will integrate
with an appropriate time-stepping scheme (detailed in section \ref{newmark}). We note that this finite
element formulation is a standard approach for linear elastodynamics and has
been widely used in biomechanical simulations and MRE studies (e.g.,
\cite{honarvar_comparison_2016}).

\subsubsection{Weak Enforcement of Dirichlet Boundary Conditions}

In our simulations, time-harmonic shear waves are induced via prescribed
oscillatory displacements on part of the boundary (Dirichlet BC). Rather than
enforcing these boundary conditions strongly (which can be challenging in
certain unfitted mesh or coupling scenarios), we impose them weakly by adding a
penalty term to the variational formulation. The penalty method augments the
weak form \eqref{weakform_final} with a term that penalizes deviation from the
desired boundary values. Specifically, for a prescribed displacement
$\mathbf{g}(\mathbf{x},t)$ on $\Gamma_D$, we write the continuity equation in
the space $V$, and add a penalization term
\begin{equation}
    \label{eq:penalty}
    \frac{\alpha_{\text{pen}}}{h} \int_{\Gamma_D} 
    [\mathbf{u} - \mathbf{g}] \cdot \mathbf{v} \, d\Gamma,
\end{equation}
where $h$ is a characteristic element size on $\Gamma_D$, and
$\alpha_{\text{pen}}$ is a penalty parameter. This term drives $\mathbf{u}$
toward $\mathbf{g}$ on the boundary by penalizing any mismatch. In practice, it
contributes additional entries to the stiffness matrix and force vector
corresponding to the boundary degrees of freedom. As $\alpha_{\text{pen}}/h \to
\infty$, one recovers the strong enforcement of $\mathbf{u}=\mathbf{g}$.
However, excessively large penalties can lead to ill-conditioned stiffness
matrices and numerical instability. Therefore, $\alpha_{\text{pen}}$ must be
chosen sufficiently high to enforce the displacement approximately, yet not so
high as to spoil the conditioning of the system (following guidelines from,
e.g., \cite{barrett_finite_1986,benzaken_constructing_2024}). In summary, the
penalty-enforced weak form is: find $\mathbf{u}^h\in V^h$ such that
\begin{equation}
    \label{weakform_penalty}
    \begin{aligned}
    &\int_{\Omega}\rho\ddot{\mathbf{u}}^h\cdot \mathbf{v}^h d\Omega +
    \int_{\Omega}2\mu\boldsymbol{\epsilon}(\mathbf{u}^h):\boldsymbol{\epsilon}(\mathbf{v}^h)d\Omega
    + \int_{\Omega}2\eta\dot{\boldsymbol{\epsilon}}(\mathbf{u}^h):
    \boldsymbol{\epsilon}(\mathbf{v}^h)d\Omega \\
    &\qquad + \frac{\alpha_{\text{pen}}}{h}\int_{\Gamma_D}[\mathbf{u}^h-\mathbf{g}]\cdot
    \mathbf{v}^h d\Gamma = 0,
    \end{aligned}
\end{equation}
for all $\mathbf{v}^h\in V^h$. The corresponding matrix system
\eqref{FEM_matrix} is modified by additional penalty terms $\mathbf{K}_{\text{pen}}$ and
$\mathbf{f}_{\text{pen}}$ (stiffness and force contributions) affecting only the
constrained nodes.

\subsubsection{Absorbing Boundary: Perfectly Matched Viscous Layer}

To avoid nonphysical reflections of outgoing shear waves at the domain
boundaries, we surround the region of interest with an absorbing layer. We
implement a Perfectly Matched Layer (PML) in the form of a viscous damping
region. In practice, this is achieved by increasing the material damping within
a thin layer of thickness $H_l$ at the model boundaries, so that waves entering
this layer are gradually attenuated. We introduce a damping scaling factor
$\alpha_{\text{abs}}$ and augment the global damping matrix as
\begin{equation}
    \label{PML_damping} \mathbf{C} \leftarrow \mathbf{C} + \alpha_{\text{abs}}\mathbf{C}_{\text{layer}},
\end{equation} 
where $\mathbf{C}_{\text{layer}}$ is the contribution of elements in the absorbing layer
(we often take $\mathbf{C}_{\text{layer}} = \frac{1}{H_l}\mathbf{C}$ for those elements,
distributing the effect over the layer thickness). In other words, the viscosity
$\eta$ is effectively increased by a factor $\alpha_{\text{abs}}/H_l$ in the
boundary layer. A larger $\alpha_{\text{abs}}$ yields stronger absorption
(greater wave attenuation) but can also cause some reflection at the interface
between the physical domain and the absorbing layer if the impedance mismatch is
too abrupt. Therefore, $\alpha_{\text{abs}}$ is tuned to balance minimal
reflection with sufficient dissipation of the wave before it reaches the outer
boundary. This approach acts as a simple sponge layer and is a common technique
to simulate an open (infinite) domain for wave propagation problems (cf. the
original PML concept by Berenger for electromagnetic waves
\cite{berenger_perfectly_1994} and its adaptations in elastodynamics
\cite{johnson_notes_2021}).

\subsubsection{Time Discretization: Newmark-Beta Method}\label{newmark}

The semi-discrete equations \eqref{FEM_matrix} are integrated in time using the
Newmark-$\beta$ method. Newmark's method is an implicit, second-order time
integration scheme widely used in structural dynamics due to its stability and
accuracy properties \cite{Nwmark_method_1959,hughes_finite_1987}. We update the
nodal displacements $\mathbf{d}$ and velocities $\dot{\mathbf{d}}$ over a time
step $\Delta t$ according to the Newmark formulas: 
\begin{equation}
\begin{aligned}
    \mathbf{d}_{n+1} &= \mathbf{d}_n + \Delta t\,\dot{\mathbf{d}}_n + \frac{\Delta t^2}{2}[1-2\beta]\,\ddot{\mathbf{d}}_n + \beta\,\Delta t^2\,\ddot{\mathbf{d}}_{n+1}~, \\
    \dot{\mathbf{d}}_{n+1} &= \dot{\mathbf{d}}_n + [1-\gamma]\,\Delta t\,\ddot{\mathbf{d}}_n + \gamma\,\Delta t\,\ddot{\mathbf{d}}_{n+1}~.
\end{aligned}
\label{newmark_disp_vel}
\end{equation}
where $\beta$ and $\gamma$ are the Newmark scheme parameters. We choose the
classical values $\gamma=\tfrac{1}{2}$ and $\beta=\tfrac{1}{4}$, corresponding
to the constant average acceleration method, which is unconditionally stable and
second-order accurate for linear systems. At each time step,
\eqref{newmark_disp_vel} is used in conjunction with the equilibrium equation
\eqref{FEM_matrix} to solve for the unknown acceleration
$\ddot{\mathbf{d}}_{n+1}$ (and hence update $\mathbf{d}_{n+1}$ and
$\dot{\mathbf{d}}_{n+1}$). In matrix form, the Newmark integration leads to a
linear solve for $\ddot{\mathbf{d}}_{n+1}$ at each step:
\begin{equation}
    \label{newmark_matrix}
    \left[ \mathbf{M} + \gamma \Delta t\, \mathbf{C} + \beta \Delta t^2\, \mathbf{K} \right] \ddot{\mathbf{d}}_{n+1} = \mathbf{f}^{\text{eff}}_{n+1},
\end{equation}
where $\mathbf{f}^{\text{eff}}_{n+1}$ is an effective load vector including
external forces at $t_{n+1}$ and contributions from known quantities at $t_n$.
This implicit scheme provides numerical damping of high-frequency modes and
ensures stability even for large time steps, which is advantageous for the
transient simulations in this study.

\subsubsection{Vascular Inclusion}\label{RLM}

To model a pulsating blood vessel embedded in the tissue, we employ a multiscale
3D-1D coupling approach based on a reduced Lagrange multiplier (RLM) method
\cite{heltai_multiscale_2019, belponer_reduced_2023, heltai_zunino_2023}. The blood vessel is
represented as a 1D line inclusion (along its centerline) within the 3D domain,
and its mechanical effect on the tissue is introduced via coupling conditions on
the vessel boundary. In essence, a time-dependent internal pressure $p(t)$ is
prescribed inside the vessel, and this exerts a Neumann (traction) boundary
condition on the surrounding tissue along the vessel wall. We impose the
constraint 
\begin{equation}\label{vessel_BC} 
    \boldsymbol{\sigma}(\mathbf{u}) \cdot \mathbf{n}_v = -p(t) \mathbf{n}_v~, 
\end{equation} 
on the vessel boundary (here $\mathbf{n}_v$ denotes the outward normal from the
vessel into the tissue). This condition is enforced weakly using Lagrange
multipliers defined along the 1D inclusion, following the formulation in
\cite{heltai_multiscale_2019}. Physically, \eqref{vessel_BC} causes the vessel
to expand and contract uniformly with the applied pressure (since the traction
is normal and proportional to $p(t)$), transmitting cyclic deformation to the 3D
mesh. By coupling the 1D pressure model with the 3D solid mechanics, we capture the essential effect of vascular pulsation in MRE: a localized, periodic prestress and displacement field that modifies shear-wave propagation in the surrounding tissue. Although the constitutive material parameters remain unchanged, this interaction leads to an apparent increase in the effective stiffness recovered by inversion algorithms, and it is reflected in the elastograms through altered wavelengths and interference patterns. This coupled modeling strategy therefore enables us to study how vascular pulsation influences MRE-based parameter reconstructions within a purely linear material framework.

    \subsubsection{Algorithm}
    The aforementioned model is implemented in C++ using deal.ii v9.5.0 \cite{arndt_dealii_2023}. It is run on the High Performance Computing Cluster, Meggie, having 728 compute nodes each with two intel ``Broadwell'' chips. The algorithm is as shown in \ref{alg:forward_prob}.
    \begin{algorithm}
    \caption{Shear Wave Propagation in Linear Isotropic Viscoelastic Medium}
    \begin{algorithmic}[1]
        \Require Input: \texttt{Mesh, Initial Conditions, Material Parameters}
        \For{each element $e = 1$ to $N_{\text{elements}}$}
            \For{each local degree of freedom $a = 1$ to $n_{\text{local}}$}
                \For{each local degree of freedom $b = 1$ to $n_{\text{local}}$}
                    \State $K_{ab}^{(e)} = \int_{\Omega_e} 2\mu \boldsymbol{\epsilon} (\boldsymbol{N}^a)\cdot \boldsymbol{\epsilon} (\boldsymbol{N}^b)d\Omega_e$
                    \State $C_{ab}^{(e)} = \int_{\Omega_e} 2\eta \boldsymbol{\epsilon} (\boldsymbol{N}^a) \cdot \boldsymbol{\epsilon} (\boldsymbol{N}^b)d\Omega_e$
                    \State $M_{ab}^{(e)} = \int_{\Omega_e} \rho \boldsymbol{N}^a \cdot \boldsymbol{N}^bd\Omega_e$
                \EndFor
                \State $f_{a}^{(e)} = \int_{\Omega_e} f\boldsymbol{N}^ad\Omega_e$
            \EndFor
        \EndFor
        \State Assemble global matrices: $\mathbf{K} = \mathbb{A}(\mathbf{K}^{(e)})$, $\mathbf{C} = \mathbb{A}(\mathbf{C}^{(e)})$, $\mathbf{M} = \mathbb{A}(\mathbf{M}^{(e)})$
        \State Apply $\boldsymbol{\sigma}(\mathbf{u})\cdot\mathbf{n}=-p\mathbf{n}$ at the domain-inclusion boundary
        \State Apply $\frac{\alpha_{\rm pen}}{h}$ penalty at the external Dirichlet boundary 
        \State Apply $\frac{\alpha_{\rm abs}}{H_l}$ viscosity at the absorbing Perfectly Matched Layer
        \State Initialize $\mathbf{u}_0 = \bold{0}$, $\dot{\mathbf{u}}_0 = \bold{0}$
        \For{each time step $t = 0$ to $T$}
            \State Predictor: $\tilde{\mathbf{u}}_{t+1} = \mathbf{u}_t + \dot{\mathbf{u}}_t \Delta t + \ddot{\mathbf{u}}_t [\frac{1}{2} - \beta] \Delta t^2$
            \State Predictor: $\dot{\tilde{\mathbf{u}}}_{t+1} = \dot{\mathbf{u}}_t + \ddot{\mathbf{u}}_t [1 - \gamma] \Delta t$
            \State Solve: $\ddot{\mathbf{u}}_{t+1} = [\mathbf{M} + \mathbf{C} \gamma \Delta t + \mathbf{K} \beta \Delta t^2]^{-1} \left[\mathbf{f}_{t+1} - \mathbf{K} \tilde{\mathbf{u}}_{t+1} - \mathbf{C} \dot{\tilde{\mathbf{u}}}_{t+1}\right]$
            \State Corrector: $\mathbf{u}_{t+1} = \tilde{\mathbf{u}}_{t+1} + \beta \ddot{\mathbf{u}}_{t+1} \Delta t^2$
            \State Corrector: $\dot{\mathbf{u}}_{t+1} = \dot{\tilde{\mathbf{u}}}_{t+1} + \gamma \ddot{\mathbf{u}}_{t+1} \Delta t$
        \EndFor
        \Require Output: \texttt{Displacement Field $\mathbf{u}(\mathbf{x}, t)$}
    \end{algorithmic}
    \label{alg:forward_prob}
    \end{algorithm}

\subsection{Direct Inversion}\label{Inversion}

Local algebraic inversion is a widespread and popular technique for reconstruction of material properties from wave data. This Direction Inversion (DI) method is based on the
isotropic wave propagation equation. In a linear viscoelastic medium, wave
propagation obeys the law of linear momentum (derived earlier in
algorithm ~\ref{alg:forward_prob} and equation ~\ref{cauchy_eq}, ~\ref{sigma_eq}). The equation of
motion for a (possibly inhomogeneous), isotropic, linear viscoelastic medium is:
\begin{equation}
\begin{aligned}
    \rho\,\ddot{\mathbf{u}}(\mathbf{x}, t) =\;&
    \nabla \cdot \left( \mu\, \nabla \mathbf{u}(\mathbf{x}, t) \right)
    + \nabla \left( (\lambda + \mu)\, \nabla \cdot \mathbf{u}(\mathbf{x}, t) \right) \\
    &+ \nabla \cdot \left( \eta\, \nabla \dot{\mathbf{u}}(\mathbf{x}, t) \right)
    + \nabla \left( \eta\, \nabla \cdot \dot{\mathbf{u}}(\mathbf{x}, t) \right),
\end{aligned}
\end{equation}
where $\rho$ is the density, $\mu$ and $\lambda$ are the Lamé parameters (shear modulus and first Lamé constant), and $\eta$ is the viscosity coefficient. Assuming $\mu$, $\lambda$  and $\eta$ are constant in space and time, the above equation simplifies to:
\begin{equation}
    \rho\,\ddot{\mathbf{u}}(\mathbf{x}, t) = \mu\,\Delta \mathbf{u}(\mathbf{x}, t) + [\lambda + \mu]\nabla(\nabla \cdot \mathbf{u}(\mathbf{x}, t)) + \eta\,\Delta \dot{\mathbf{u}}(\mathbf{x}, t) + \eta\,\nabla(\nabla \cdot \dot{\mathbf{u}}(\mathbf{x}, t)),
\end{equation}
which still includes both shear (transverse) and compressional (longitudinal) wave components. If no longitudinal waves are present (i.e. the displacement field is divergence-free, $\nabla\cdot \mathbf{u}=0$), the compressional terms vanish. Enforcing this Helmholtz decomposition (to ensure a divergence-free $\mathbf{u}_{\text{curl}}$) further reduces the governing equation to:
\begin{equation}
    \rho\,\ddot{\mathbf{u}}_{\text{curl}}(\mathbf{x},t) = \mu\,\Delta \mathbf{u}_{\text{curl}}(\mathbf{x},t) + \eta\,\Delta \dot{\mathbf{u}}_{\text{curl}}(\mathbf{x},t),
\end{equation}
which describes purely transverse wave propagation. Taking a temporal Fourier
transform of this equation and isolating a single steady-state excitation
frequency $\omega$ yields a frequency-domain form. In the frequency domain, time
derivatives translate to factors of $i\omega$, and we obtain a Helmholtz-type
equation for the complex vector valued field
$\tilde{\mathbf{u}}_{\text{curl}}(\mathbf{x}, \omega)$:
\begin{equation}
    \label{helmholtz_eq}
    \rho \omega^2 \tilde{\mathbf{u}}_{\text{curl}}(\mathbf{x}, \omega) + [\mu + i \omega \eta] \Delta \tilde{\mathbf{u}}_{\text{curl}}(\mathbf{x}, \omega) = 0,
\end{equation}
where $\tilde{\mathbf{u}}(x,\omega)$ denotes the complex-valued
Fourier-transformed displacement field at angular frequency $\omega$. 

Equation~\eqref{helmholtz_eq} is valid for homogeneous media, where the shear
modulus $\mu$ and viscosity $\eta$ are constant throughout the domain. In case
of heterogeneous media, one can still use the same equation, provided that $\mu$
and $\eta$ are constants within each local region of interest. In this case,
such equation is still valid away from material interfaces and boundaries, i.e.,
where $\partial \eta/\partial\mathbf{x} = 0$ and $\partial \mu/\partial \mathbf{x} = 0$.

In these regions, one may rearrange \eqref{helmholtz_eq} to obtain an
overdetermined direct relationship between the time-Fourier transform of the
measured displacement field $\tilde{\mathbf{u}}_{\text{curl}}(\mathbf{x},
\omega)$ and the complex shear modulus 
\begin{equation}
    G^*(\omega) := \mu + i\omega\eta, 
\end{equation}
i.e., at every point $\mathbf{x}$ away from interfaces and boundaries, we satisfy the six (three each for the real and imaginary parts)
linear equations:
\begin{equation}
    \label{direct_inversion_eq}
    \begin{aligned}
        &G^*(\omega)= \frac{ \rho\,\omega^2\,[\mathbf{\tilde{u}_\text{curl}}(\mathbf{x}, \omega)]_j}{ [\Delta\mathbf{\tilde{u}}_\text{curl}(\mathbf{x}, \omega)]_j}, & j = 1,2,3.
    \end{aligned}
\end{equation}

This expression provides a direct formula for the complex shear
modulus at each location. In practice, the implementation of this relationship
can be challenging due to noise and other artifacts in the measured displacement
field, and the local inversion in ~\eqref{direct_inversion_eq} is averaged in
each coordinate direction, and in each point belonging to a subregion to
mitigate noise and the ill-posedness of the inverse problem. Under the
assumption of piecewise constant local homogeneity, this relationship can be
used to directly compute the viscoelastic material parameters from the known
displacement field at each location away from boundaries and interfaces. In a Kelvin-Voigt material model, the real part (storage
modulus) is given by $G' = \mu$  and the imaginary part (loss modulus) by $G'' = \omega\eta$. In
practice, implementing this algebraic DI requires calculating the
Laplacian of the measured displacement field. A straightforward approach is to
use finite differences. For example, the standard FD stencil for the second derivative in the $x$ coordinate direction reads as
\begin{equation}\label{laplacian}
    \partial_{xx} \tilde{\mathbf{u}}_{\text{curl}}(\mathbf{x},\omega) \approx \frac{\tilde{\mathbf{u}}_{\text{curl}}(\mathbf{x}+h\mathbf{e}_x,\omega) - 2\tilde{\mathbf{u}}_{\text{curl}}(\mathbf{x},\omega) + \tilde{\mathbf{u}}_{\text{curl}}(\mathbf{x}-h\mathbf{e}_x,\omega)}{h^2},
\end{equation}
where $h$ is the grid step size in the given coordinate direction, and $\mathbf{e}_x$ is the unit vector along $x$. The full Laplacian is then $\Delta \tilde{\mathbf{u}}_{\text{curl}} = \partial_{xx}\tilde{\mathbf{u}}_{\text{curl}} + \partial_{yy}\tilde{\mathbf{u}}_{\text{curl}} + \partial_{zz}\tilde{\mathbf{u}}_{\text{curl}}$, applied component-wise. For the concrete implementation of the FD Laplacian in three dimensions, in this study we follow Papazoglou et al.\cite{papazoglou_algebraic_2008} by employing the gradient function provided by MATLAB recursively. This function evaluates the central difference for interior data points as
\begin{equation}\label{grad_matlab}
    \partial_x \tilde{\mathbf{u}}_{\text{curl}}(\mathbf{x},\omega) \approx \frac{\tilde{\mathbf{u}}_{\text{curl}}(\mathbf{x}+h\mathbf{e}_x,\omega) -  \tilde{\mathbf{u}}_{\text{curl}}(\mathbf{x}-h\mathbf{e}_x,\omega)}{2h},
\end{equation}
and similarly for the $y$ and $z$ directions.

For the second derivative nested gradient evaluation for the interior data points gives an approximation of 
\begin{equation}\label{lap_matlab}
    \partial_{xx} \tilde{\mathbf{u}}_{\text{curl}}(\mathbf{x},\omega) \approx \frac{\tilde{\mathbf{u}}_{\text{curl}}(\mathbf{x}+2h\mathbf{e}_x,\omega) -  2\tilde{\mathbf{u}}_{\text{curl}}(\mathbf{x},\omega) + \tilde{\mathbf{u}}_{\text{curl}}(\mathbf{x}-2h\mathbf{e}_x,\omega)}{4h^2}.
\end{equation}
Observe that the resulting stencil differs from the standard FD stencil in equation \ref{laplacian}, and the full Laplacian is obtained by summing the corresponding second derivatives in $x$, $y$, and $z$.
In the current study, we apply the stencil in equation equation \ref{lap_matlab} and assume $h$ to be the size of the element in the FEM mesh, implying that the resolution of the MRE inversion is the same as the mesh resolution.

\section{Results}
\begin{figure}[h]
        \includegraphics[width=\textwidth]{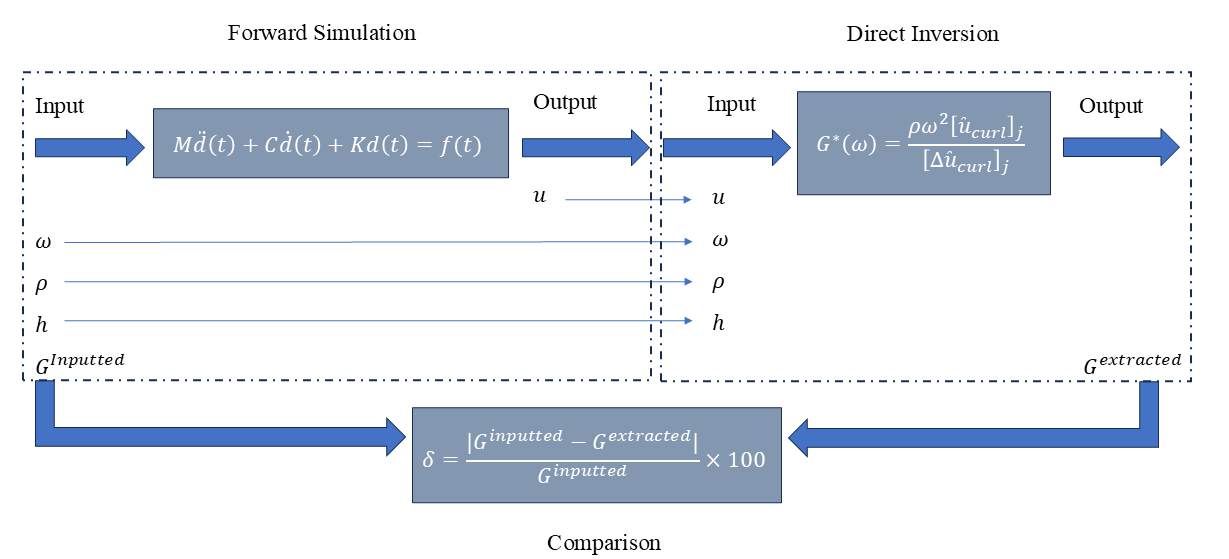}
        \centering
        \caption{Validation scheme for benchmark simulations }
        \label{fig:loop}
    \end{figure}

Various factors play a role in the successful reconstruction of material parameters using MRE. This study aims to further this understanding by using synthetic data from FE modeling. A model simulating a shear wave in a 3D cube of known material properties is developed as detailed in section \ref{FEmodel}. The displacement field from the model is inputted into the DI scheme (detailed in section \ref{Inversion}) and the re-extracted material parameters are compared to those inputted. This forms a closed-loop validation and allows us to demonstrate use of the proposed benchmarking scheme. This validation loop is depicted in figure \ref{fig:loop}.

This allows for a controlled isolation of individual effects on shear modulus reconstruction using DI. The following subsections focus on the effects of spatial and temporal discretization, material domain heterogeneity, vascular pulsation and a combination of such.

    \subsection{Spatial and Temporal Mesh Convergence}
\begin{figure}[H]
        \includegraphics[width=\textwidth]{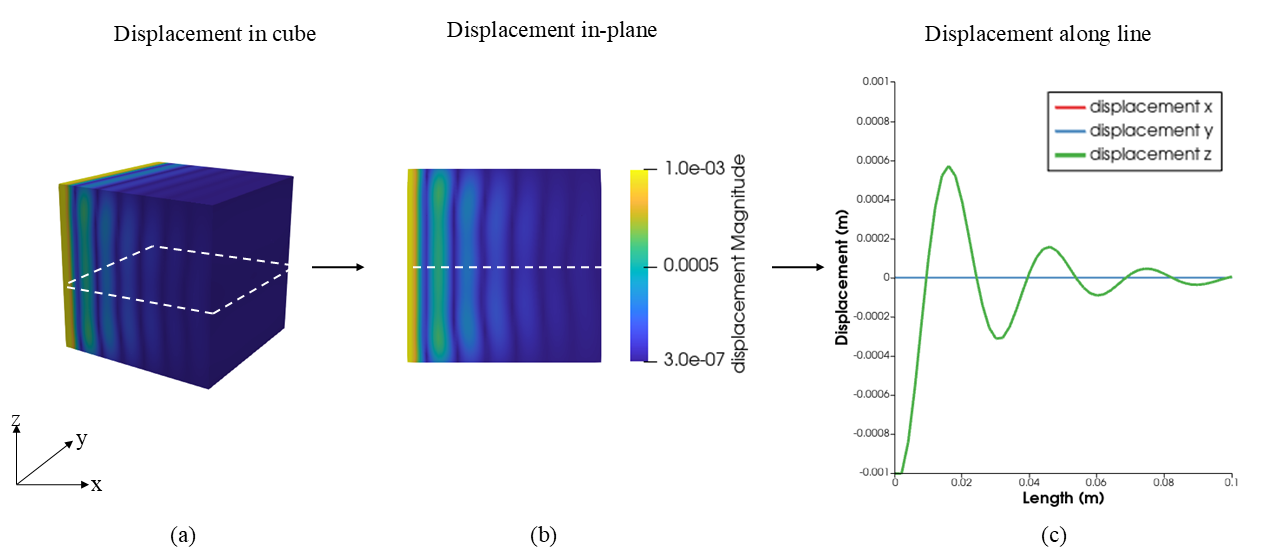}
        \centering
        \caption{(a) Magnitude of Displacement Field in an isotropic cuboidal domain, (b) Magnitude of Displacement Field in the plane highlighted in (a), (c) Displacement Field along the line highlighted in (b)  }
        \label{fig:displacement_field}
    \end{figure}

    A cube of 0.1 m x 0.1 m x 0.1 m is uniformly meshed and a shear wave of 50 Hz is applied on  one boundary to emulate MRE loading conditions. The remaining boundaries have an absorbing layer to minimize the reflection in the system. The displacement field so obtained is depicted in figure \ref{fig:displacement_field}.

\begin{figure}[h]
      \begin{subfigure}[b]{0.5\textwidth}
        \includegraphics[width=\textwidth]{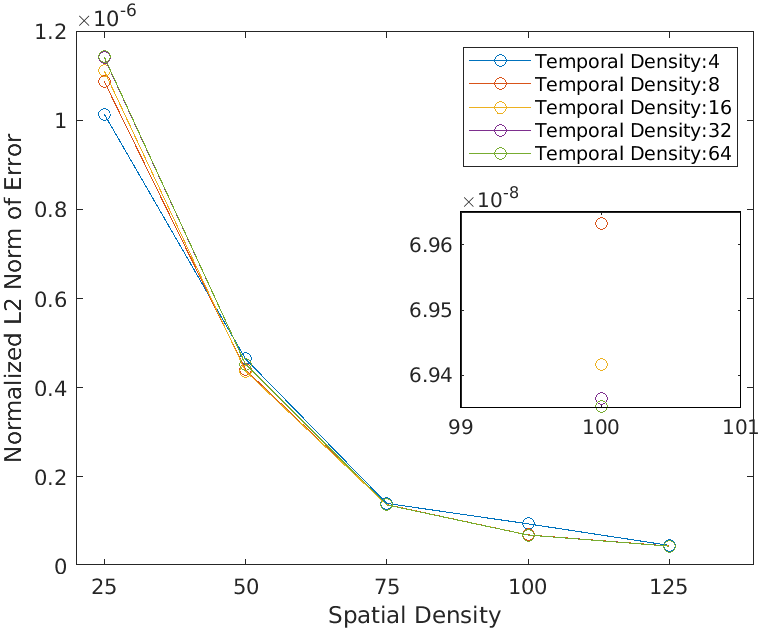}
        \caption{}
        \label{fig:spatial_conv}
      \end{subfigure}
      \hfill
      \begin{subfigure}[b]{0.5\textwidth}
        \includegraphics[width=\textwidth]{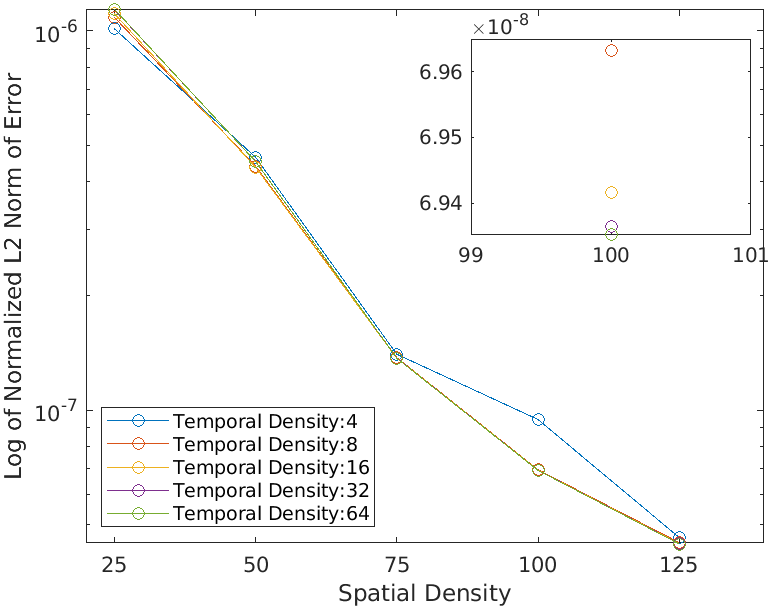}
        \caption{}
        \label{fig:spatial_conv_log}
      \end{subfigure}
      \caption{Convergence of error (\ref{fig:spatial_conv}) and logarithmic convergence of error (\ref{fig:spatial_conv_log}) for refining the mesh resolution at different temporal densities}
      \label{spatial convergence}
    \end{figure}

\begin{figure}[!tbp]
      \begin{subfigure}[b]{0.5\textwidth}
        \includegraphics[width=\textwidth]{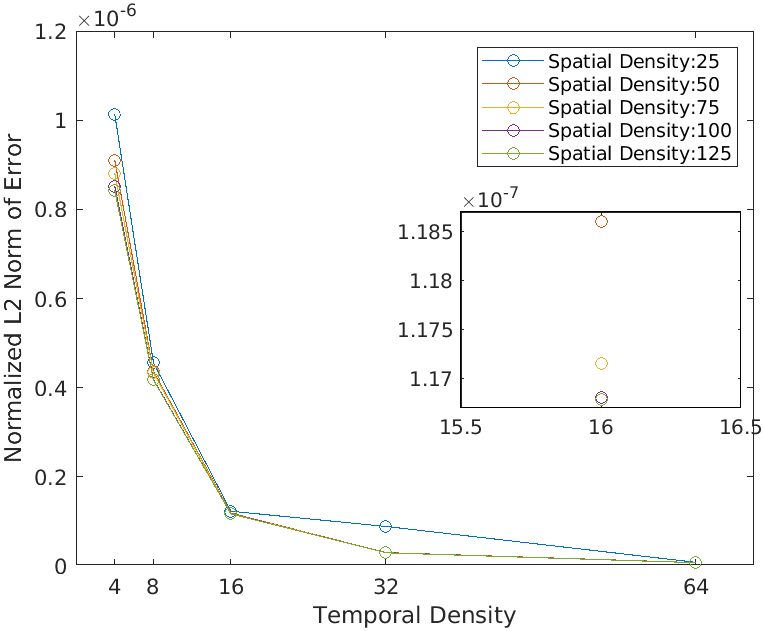}
        \caption{}
        \label{fig:temporal_conv}
      \end{subfigure}
      \hfill
      \begin{subfigure}[b]{0.5\textwidth}
        \includegraphics[width=\textwidth]{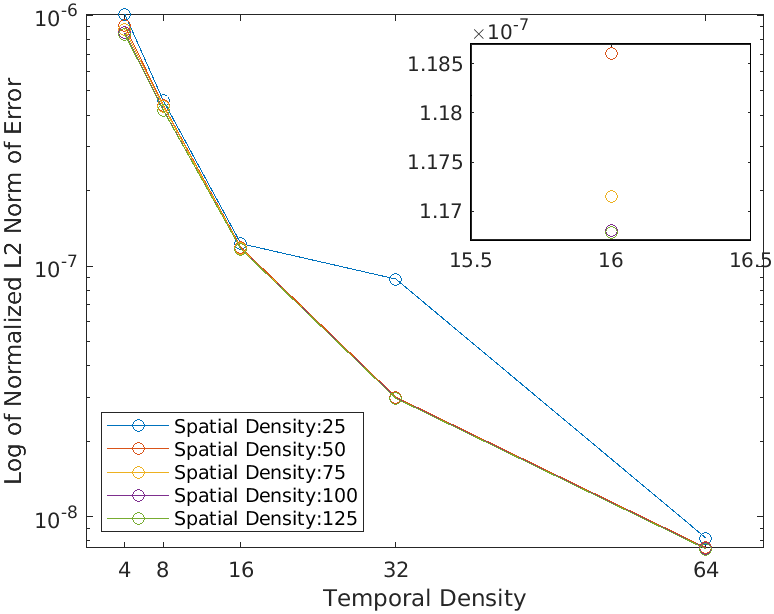}
        \caption{}
        \label{fig:temporal_conv_log}
      \end{subfigure}
      \caption{Convergence of error (\ref{fig:temporal_conv}) and logarithmic convergence of error (\ref{fig:temporal_conv_log}) for refining the temporal resolution at different mesh densities}
      \label{temporal convergence}
    \end{figure}

    Since this forward elastodynamic problem is discretized in space and time, first a mesh convergence study is performed to ensure the reliability and accuracy of the simulated results. Firstly, for each fixed temporal density, the mesh is systematically refined in-plane. This is done by incrementally increasing the number of nodes in each direction from 25 to 125, in steps of 25. This increases the spatial resolution from 7.325 to 14.65, 22.538, 29.3 and finally 36.635 pixels per wavelength. Secondly, for each spatial grid, the sampling time is increased from $2^2$ to $2^6$ in 4 increments. At each refinement level, the displacement at the nodes of mesh $\mathbf{u}^{N-1}$ are compared to the displacements of the same nodes in a denser mesh $\mathbf{u}^{N}$ using the L2 norm. The difference $||\triangle \mathbf{u}^{N}||_0=||\mathbf{u}^{N}-\mathbf{u}^{N-1}||_0$ is normalized and plotted as shown in figure \ref{spatial convergence} and \ref{temporal convergence}. In the forward wave simulation, error decreases with refinement.
    
\begin{figure}[H]
        \centering
        \begin{subfigure}[b]{0.49\textwidth}
            \centering
            \includegraphics[width=\linewidth]{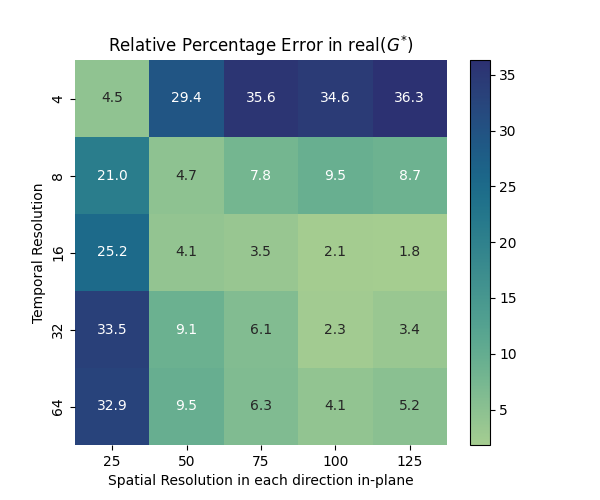}
            \caption{}
            \label{fig:FD_real}
        \end{subfigure}
        \hfill
        \begin{subfigure}[b]{0.49\textwidth}
            \centering
            \includegraphics[width=\linewidth]{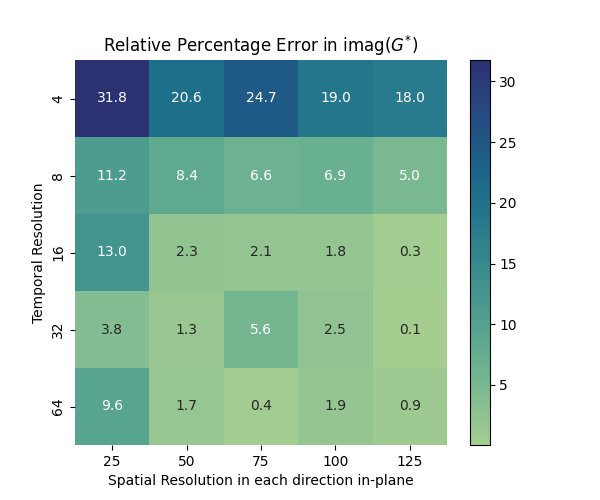}
            \caption{}
            \label{fig:FD_imag}
        \end{subfigure}
        \vspace{0.1cm} 
        \caption{Relative Percentage error in the Storage (\ref{fig:FD_real}) and loss modulus (\ref{fig:FD_imag}) for the DI method at different temporal and mesh refinements}
        \label{fig:conergence heatmaps}
    \end{figure}

     However, the DI reconstruction shows non-monotonic convergence, indicating that finer resolutions do not always reduce inversion error. The reconstructed $G^{extracted}$ is averaged over the mesh and in doing so $G^{extracted}$ values at the boundary nodes are disregarded to avoid inaccuracies. The recovered material parameters are compared to those inputted by calculating a relative percentage error as 
    \begin{equation}
    \delta =\frac{|G^{inputted}-G^{extracted}|}{G^{inputted}}\times100.
    \end{equation} 
    This is shown in figure  \ref{fig:conergence heatmaps}. The non-monotonic convergence is attributed to two sources of numerical errors when calculating the Laplacian in the complex plane using the central difference approximation, namely, truncation and subtractive cancellation \cite{javili_computational_2020}. These dominate at different mesh resolutions \cite{lai_extensions_2008},\cite{fike_optimization_2011}. 
Thus, there exists an ideal mesh based on the material parameters and frequency of external wave. Hence, adjusting the spatial and temporal resolution of the MRE, ipso facto the voxel size of the measurement based on the wavelength of the shear displacement field can provide accurate results with lower computational costs. Moreover, the sign of this error is of additional consideration, since it can determine the dominant error sources as noise at high spatial support leads to systematic underestimation of shear moduli whereas coarse spatial discretization induces overestimation \cite{mura_analytical_2020}.
      
\begin{figure}[!tbp]
        \centering
        \includegraphics[width=0.5\textwidth]{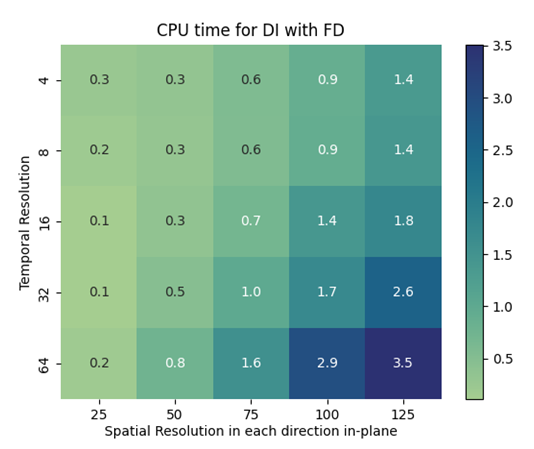}
        \caption{Computational cost of the DI method in seconds for different temporal and spatial refinements}
        \label{CPU hrs}
    \end{figure}

    Figure \ref{CPU hrs} shows the computational cost in CPU time for the inversion. Mesh refinement increases computational time more significantly, due to point-wise computation of the Laplacian, whereas temporal refinement primarily adds an overhead when taking the Fourier transform.
    
    \subsection{Material Domain Heterogeneity}\label{hetero}
\begin{figure}[h]
        \includegraphics[width=\textwidth]{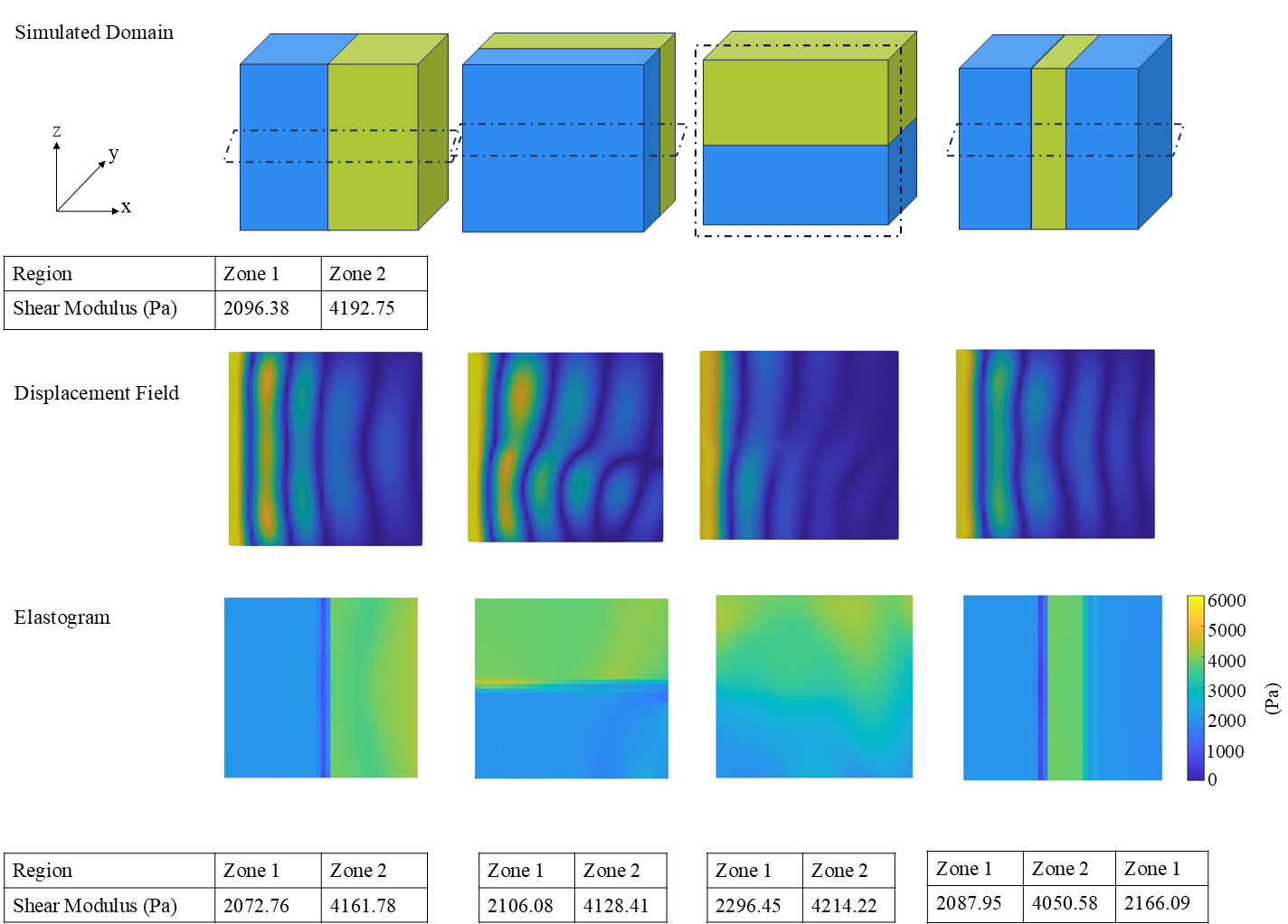}
        \centering
        \caption{Simulated domains with multiple zones and the respective magnitudes of the displacement fields and elastograms in the highlighted planes}
        \label{fig:zones}
    \end{figure}

    For the next study, the cuboidal domain is divided into different zones and each zone is given a different material property as can be seen in figure \ref{fig:zones}. The shear modulus is kept constant within each zone such that the heterogeneous domain is actually a composition of homogeneous zones. The displacement field from the forward simulation are fed into the inversion algorithm and the so developed elastograms in the highlighted planes are also shown in the same figure. In calculating the recovered shear modulus, we disregard the interface boundary of the zones, again due to inaccuracies. \\
    
    As can be seen in the figure, DI successfully reconstructed the shear moduli in the distinct zones with a maximum of 8.45 \% relative error. The loss moduli error was consistently higher than the error in storage moduli, suggesting that low-amplitude reflected waves at the zonal interface affect the viscous parameter more than the elastic one. \\
    
    Another notable aspect is the distinctiveness of the interface boundary in the elastograms. In this study, the external excitation travels along the x-axis with its amplitude oriented in the z direction. Consequently, the changes in the material properties along the z-axis are not sharply visible. The drop in the shear modulus at the interface is effected by the sharp changes in gradient when the wave directly strikes the interface. As a result, the material boundary becomes more pronounced when the wave interacts with it perpendicularly.\\

    \subsection{Pulsating Inclusions}

\begin{figure}[h]
        \includegraphics[width=\textwidth]{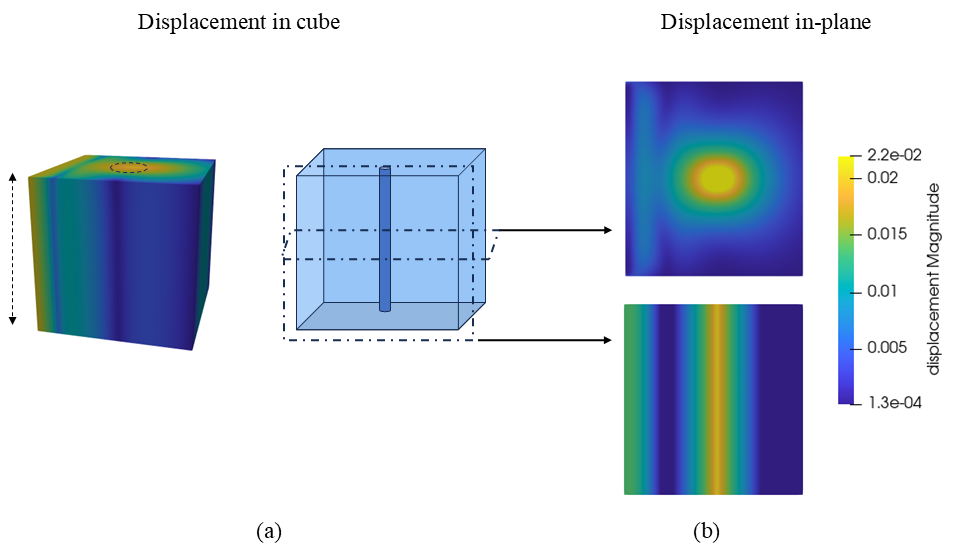}
        \centering
        \caption{Simulated domains with inclusions in different directions and respective elastograms }
        \label{fig:vessel_displacement}
    \end{figure}

    For this study, an inclusion with specified pressure and sinusoidal pulsation in time is introduced and coupled with the cuboidal domain in the aforementioned way. This inclusion is allowed to expand and contract under the excitation pressure. This introduces two interacting wave fields into the system, i.e., the external transverse wave enforced to apply the MRE testing modality and the internal wave mimicking an \textit{in-vivo} vascular pulsation. The displacement fields so produced are depicted in figure \ref{fig:vessel_displacement} in two different planes.\\
    
    Notably, the pulsation makes the inclusion overtly and sharply visible in the elastogram. The shear modulus inside the body of the inclusion drops to much lower values contrasting with a general stiffening of the surrounding tissue domain. \\

    Keeping the orientation of the vessel constant, the pressure applied to the inclusion is varied in mean (12500 - 17500 Pa) and amplitude ($\pm 2000$- $\pm 3000$ Pa). The observed material changes are listed in table \ref{tab:pressure_vessel}. Both mean and amplitude increase, raises the observed stiffness, with mean pressure producing a stronger effect. While higher amplitudes increase the force gradients, a higher mean pressure adds more overall force to the system.\\

    \begin{table}[h]
        \centering
        \small
        \begin{minipage}{0.49\textwidth}
            \centering
            \begin{tabular}{|m{9em}|m{9em}|}
                \toprule
                 Pressure in Inclusion (Pa)& Shear Modulus (Pa)\\
                 \midrule
                 12500 $\pm 2000$ & 2506.64 \\
                 15000 $\pm 2000$ & 2682.58 \\
                 17500 $\pm 2000$ & 2721.92 \\
                 \bottomrule
            \end{tabular}
            
        \end{minipage}
        \hfill
        \begin{minipage}{0.5\textwidth}
            \centering
            \small
            \begin{tabular}{|m{9em}|m{9em}|}
                \toprule
                 Pressure in Inclusion (Pa) & Shear Modulus (Pa) \\
                 \midrule
                 12500 $\pm 2000$ & 2506.64 \\
                 12500 $\pm 2500$ & 2522.98\\
                 12500 $\pm 3000$ & 2539.87\\
                 \bottomrule
            \end{tabular}
            
        \end{minipage}
        \caption{Shear modulus for different pressure amplitude and mean in the inclusion}
        \label{tab:pressure_vessel}
    \end{table}

\begin{figure}[H]
        \includegraphics[width=\textwidth]{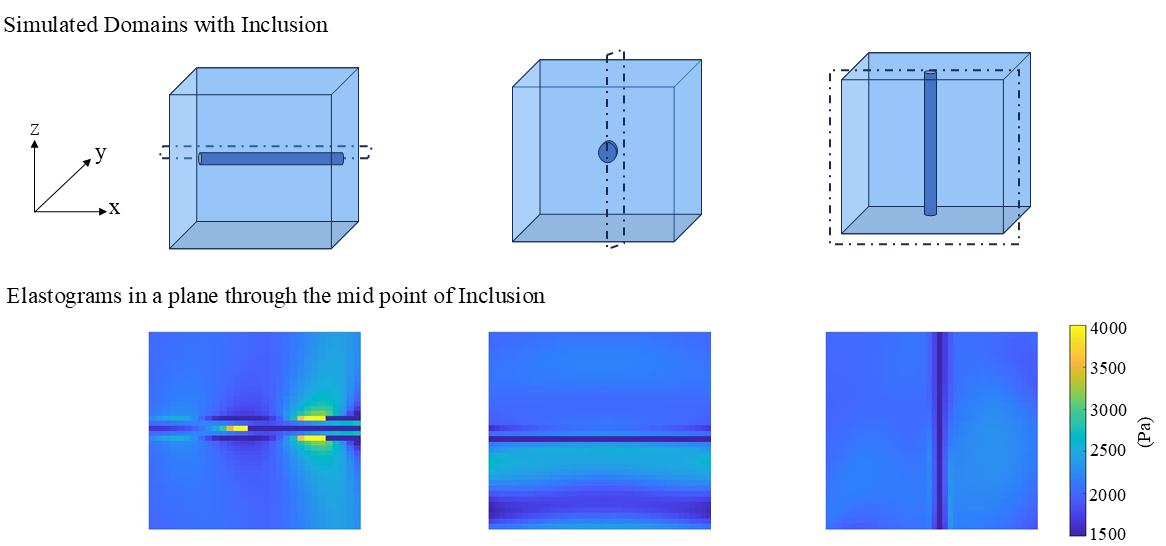}
        \centering
        \caption{Simulated domains with inclusions in different directions and respective elastograms }
        \label{fig:vessel_direction}
    \end{figure}

    Then, the orientation of the vessel with respect to the external wave was varied. This is depicted in figure \ref{fig:vessel_direction} along with the elastograms in the highlighted planes. The deformation of the surrounding tissue due to the internal pulsation interferes with the deformation due to the external excitation. Thus, depending on their relative direction, the inclusion patterns vary in the elastograms. This directional nature has implications on the shear-wave filtering steps in MRE, where directional attenuation may suppress or amplify vascular effects.\\

    The number of inclusions also have an effect on the stiffening. An increase in the number of inclusions from 1 to 2 to 4 affects the associated shear modulus that increases from 2489.27 Pa to 2580.14 Pa to 2648.48 Pa. This suggests a positive correlation between vascularization and the associated effective stiffness.\\

    \subsection{Combined Effect of Material Heterogeneity and Pulsating Inclusions}
\begin{figure}[h]
            \includegraphics[width=\textwidth]{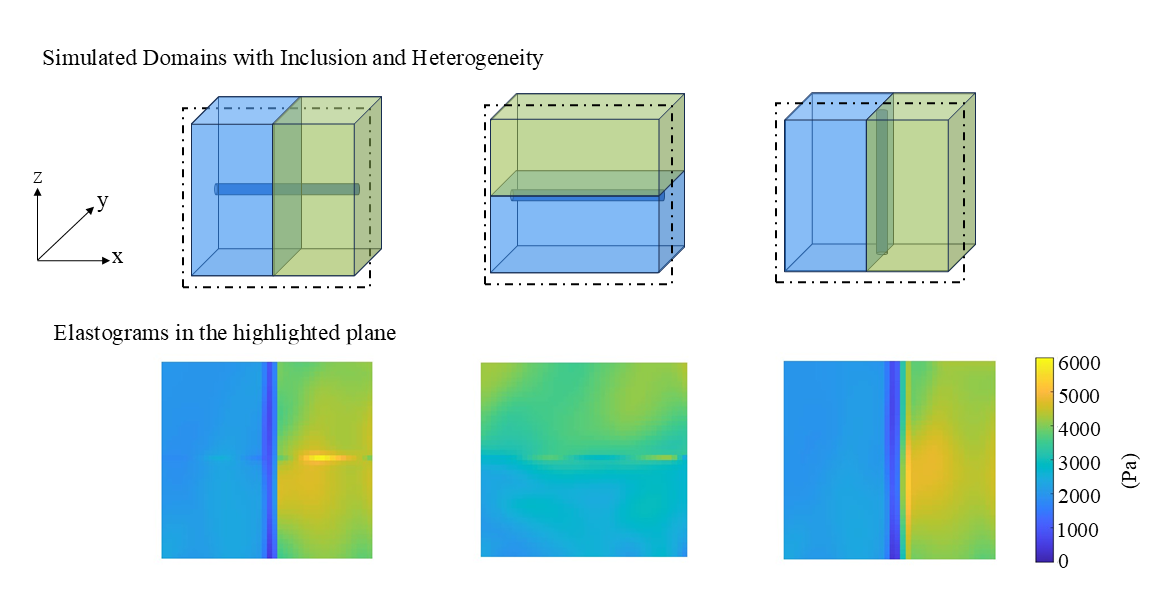}
            \centering
            \caption{Simulated domains with two zones and pulsating inclusion along  the respective elastograms}
            \label{fig:vessel_plus_hetero}
    \end{figure}

    Finally, we investigate the combined effect of both multiple material properties and vascular inclusions. The forward problem is solved on a cuboidal domain with heterogeneous domains and pulsating inclusions. The resulting displacements are used as an input for DI. The generated elastograms and the simulated domains are shown in figure \ref{fig:vessel_plus_hetero}.\\
    
    In the first case, when the vessel is oriented perpendicularly to the zonal interface, both the vessel and the interface boundary are distinctly identifiable. However, the two effects differ quantitatively, the sharp gradient changes due to change in material properties cause the shear modulus to dip sharply. While, inside the vessel, we interpolate the surrounding domain and so the gradient of change is driven solely by the additional pressure. This pressure-induced change is less intense than the change at the boundary, making the vessel appear weaker compared to the domain interface. This is accompanied with a pressure-pulse driven stiffening of the two zones. The expected shear modulus for the zones are 2096.37 Pa and 4192.74 Pa and observed to be 2192.54 Pa and 4451.24 Pa which are higher even after taking into account the relative error in our closed-loop validation in section \ref{hetero}.\\
    
    In the second case, when the vessel is oriented flush along the material interface, it coincides with the location of boundary interface, however only the vessel is visible in the elastogram. There is again a stiffening to 2511.60 Pa and 4432.81 Pa in the respective zones.\\
    
    In the final case, the vessel and zonal boundary intersect in a case where both the vessel and the interface would be overtly visible individually. However, in comparison to figure \ref{fig:zones} or the first figure in \ref{fig:vessel_plus_hetero}, the inclusion enhanced interface here is much more pronounced. This is the combined effect of the change in gradients due to the pressure application and change in material parameters. In addition, there is a stiffening of the surrounding regions. In zone 1 the shear modulus is 2356.64 Pa and in zone 2, 5506.25 Pa. Again highlighting the effect of pressure in stiffer material characterization of vascularized tissues.\\
    
\section{Discussion}
 MRE allows us to mechanically characterize soft tissues \textit{in-vivo}. The technique is based on introducing a transverse wave into the system and mapping the so-generated displacements to invert the wave equation to get the mechanical properties of the system in the form of shear moduli. This inversion of the wave equation can be performed directly, as established in section \ref{Inversion}. Such an operation allows us to extract the properties, and later associate them to known material models. However, for real materials, the validation of inversion techniques is difficult, as the material properties are a priori unknown. \\
With this study, we developed an \textit{in-silico} framework for the generation of a rich and accessible dataset that can be used as a uniform baseline for assessment of post-processing techniques in MRE based on the Finite Element (FE) method. The simulations ranged from basic homogeneous isotropic domains of varying resolution to ones with a variety of \textit{in-vivo} experimental conditions, such as heterogeneity in material composition to pulsating inclusions. The proposed validation scheme was demonstrated by comprehensive application to a DI technique.\\
The system was modeled to have linear viscoelastic properties in accordance with the Kelvin Voigt material model. This implies that the shear modulus $(G^*)$ can be decoupled into the real part that confers to the coefficient of elasticity $real(G^*)=2\mu$ and the imaginary part that correlates to the viscosity parameter as $imag(G^*)=2\pi\omega\eta$. The model was initialized with material parameters and subjected to a transverse wave, producing a displacement field that serves as input for the MRE inversion process to retrieve the material properties. By comparing the inputted and extracted material parameters, a closed-loop validation was reached for the process. In this regard, the material properties can be seen to be easily re-captured. In an attempt to reduce errors, the spatial and temporal resolution of the FE model was enhanced, but due to the nature of finite differences in complex planes, the improvement in results for better resolution was not guaranteed. This suggests that for every frequency and material parameter combination, there exists a resolution that would ensure least relative errors. Such an analysis may provide an optimum resolution for least error-prone measurements. This optimization should also consider the dependence on the signal-to-noise ratio \cite{mura_analytical_2020}\\
Soft tissues are inherently heterogeneous, and within the framework of this study, we evaluated the performance of the DI scheme in the presence of such heterogeneity. The modeled domain was split into 2 to 3 zones with different material properties. In the similar aforementioned closed-loop validation, the inputted shear moduli could be re-extracted successfully. An observation of note is to see that the boundary interfaces are very well defined in the developed elastograms, when the interface is perpendicular to the direction of travel of the transverse wave. This is due to the fact that when the wave strikes the interface, where there is a sharp change in the displacements and wave speed, the resultant gradient of change is visible in the elastrograms as a boundary artefact in the form of a dip in shear modulus values. The values lower due to the increase in the second derivative that acts as a denominator in the calculation of the same. This shear wave scattering allows detection of even minute artefacts and has been experimentally reported \cite{sack_magnetic_2023}. In the case when the heterogeneity interface is parallel to the direction of the wave, the elastograms clearly show two different zones but the boundary is not overt. This can be attributed to the fact that there is no sharp change of wave speeds, but two different wave speeds in different regions. This further motivates the need for cluster-driven algorithms to ascertain the appropriate characterization of heterogeneous regions in tissues, as the interface of internal boundaries may not be obvious. \\
 Additionally, in soft tissues, especially \textit{in-vivo}, there are a number of internal forces and at organ level, the blood pressure is one of significance. The healthy human body has a blood pressure of 120/80 mm Hg that pulses roughly at heart rate of 1 Hz. With our current study, a foundational paradigm was established to enhance understanding of the pressure-induced changes in soft tissue stiffness. This vascular pressure was exerted internally into the system by coupling the isotropic domain with inclusions using the reduced Lagrange multiplier technique in the 3D-1D domain. This implies that the inclusion is coupled such that the vessel is allowed to expand freely in response to the applied pressure, with the resulting deformation transmitted to the surrounding tissue domain, thereby influencing the displacement field. The material properties of these vessels were assumed to be the same as the surrounding domain and hence did not add any material heterogeneity by themselves. In general, the introduction of internal pressure does not alter the underlying material properties, which remain linear and unchanged, but it modifies the mechanical response of the system through additional, time-dependent internal stresses. As a consequence, the shear-wave field is altered by the superposition and interaction of externally induced and pressure-driven deformations. When processed through direct inversion, this modified wave field leads to higher apparent shear moduli compared to the input values. This behavior is therefore expected and reflects an effective stiffening inferred by the inversion procedure, rather than a physical increase in the material stiffness, arising from the presence of vascular prestress and wave interference effects.
 
 Irrespective of the direction of the vessels, their presence in the elastograms was easily ascertained. They were distinguishable with a dip in shear modulus values at the centre line. The dip was not as high as in the case of material boundaries. Increasing the pressure within the inclusion, either in mean value or amplitude, enhances the internal stress field and its interaction with the externally induced shear waves, resulting in a more pronounced effective stiffening observed in the inversion results. This effect was also enhanced based on the number of vessels, since increasing the vessels also increases the pressure. Another aspect to note is that the vessel direction can affect the elastograms differently. When the vessel is along the direction of the external transverse wave, the displacement field of the vessel adds to external wave field, this generates a much more pronounced interference pattern and hence visible vessel. In the case that the vessel is perpendicular, the elastrograms depict them as more localized lines, because the displacement field interference is perpendicular in the transverse and pressure wave. This motivates the use of filtering algorithms to remove the pressure wave effects. As the pressure is sinusoidally time-variant, the increase in the system and thereby shear moduli increase is not constant, it also has sinusoidal nature, suggesting that once the waves reach steady state, the wave pattern in time is a high frequency wave that itself oscillates in a low frequency pattern. This is also observed in the simulations.    \\
As a last step, the combined effect of the vascular inclusion and heterogeneity was taken into account in the same system. As expected and hypothesized, we see a stiffening of the zones due to added pressure. The vessel and boundary domain are clearly visible separately but when overlapped, they are visible more overtly than their individual contributions. Hence motivating the need to distinguish between interfaces and pressurized vascular sources.\\

\section{Conclusion}
A convergent and validated system was developed that can be used to \textit{in-silico} perform MRE experiments. This framework accommodates both homogeneous and heterogeneous configurations in an attempt to present a controlled platform for identification of use-case-specific inversion algorithms and accompanying raw data filtration. When employing DI schemes to heterogeneous domains, one must take into account the changes at the interface and the associated drop in observed shear moduli. In addition, the presence of vascular inclusions and associated pressure pulse gives an overall higher average stiffness in comparison to that of the surrounding tissue. Thereby indicating the need for robust filtering algorithms to remove the effect of pressurized internal waves when characterizing soft tissues.

\section*{Acknowledgments}
Funded by the German Research Foundation (DFG) project 460333672 CRC1540 EBM. LH acknowledges the MIUR Excellence Department Project awarded to the Department of Mathematics, University of Pisa, CUP I57G22000700001, and partial support from grant MUR PRIN 2022 No. 2022WKWZA8 “Immersed methods for multiscale and multiphysics problems (IMMEDIATE)”. LH is member of Gruppo Nazionale per il Calcolo Scientifico (GNCS) of Istituto Nazionale di Alta Matematica (INdAM).

\section*{Data Avalability} 
The simulation data is publicly accessible on the Zenodo repository on \url{https://doi.org/10.5281/zenodo.19236558}. The MRE processing platform is available at \url{https://bioqic-apps.charite.de/}.

\section*{Ethics Statement}
None

\bibliography{references}

@article{honarvar_comparison_2016,
    title = {A comparison of direct and iterative finite element inversion techniques in dynamic elastography},
    volume = {61},
    issn = {0031-9155, 1361-6560},
    url = {https://iopscience.iop.org/article/10.1088/0031-9155/61/8/3026},
    doi = {10.1088/0031-9155/61/8/3026},
    language = {en},
    number = {8},
    urldate = {2025-02-25},
    journal = {Physics in Medicine and Biology},
    author = {Honarvar, M and Rohling, R and Salcudean, S E},
    month = apr,
    year = {2016},
    pages = {3026--3048},
}

@article{honarvar_comparison_2017,
    title = {A {Comparison} of {Finite} {Element}-{Based} {Inversion} {Algorithms}, {Local} {Frequency} {Estimation}, and {Direct} {Inversion} {Approach} {Used} in {MRE}},
    volume = {36},
    issn = {1558-254X},
    url = {https://ieeexplore.ieee.org/document/7885109/?arnumber=7885109},
    doi = {10.1109/TMI.2017.2686388},
    number = {8},
    urldate = {2025-02-25},
    journal = {IEEE Transactions on Medical Imaging},
    author = {Honarvar, Mohammad and Sahebjavaher, Ramin S. and Rohling, Robert and Salcudean, Septimiu E.},
    month = aug,
    year = {2017},
    note = {Conference Name: IEEE Transactions on Medical Imaging},
    pages = {1686--1698},
}

@article{hollis_finite_2016,
    title = {Finite {Element} {Analysis} to {Compare} the {Accuracy} of the {Direct} and {MDEV} {Inversion} {Algorithms} in {MR} {Elastography}},
    volume = {43},
    issn = {1819-656X},
    url = {https://www.research.ed.ac.uk/en/publications/finite-element-analysis-to-compare-the-accuracy-of-the-direct-and},
    language = {English},
    number = {2},
    urldate = {2025-02-27},
    journal = {IAENG International Journal of Computer Science},
    author = {Hollis, Lyam M. and Barnhill, Eric and Conlisk, Noel and Thomas-Seale, L. E. J. and Roberts, John and Pankaj, Pankaj and Hoskins, Peter},
    month = may,
    year = {2016},
    pages = {137--146},
}

@article{barnhill_heterogeneous_2018,
    title = {Heterogeneous {Multifrequency} {Direct} {Inversion} ({HMDI}) for magnetic resonance elastography with application to a clinical brain exam},
    volume = {46},
    issn = {1361-8415},
    url = {https://www.sciencedirect.com/science/article/pii/S1361841518300835},
    doi = {10.1016/j.media.2018.03.003},
    urldate = {2025-02-26},
    journal = {Medical Image Analysis},
    author = {Barnhill, Eric and Davies, Penny J. and Ariyurek, Cemre and Fehlner, Andreas and Braun, Jürgen and Sack, Ingolf},
    month = may,
    year = {2018},
    pages = {180--188},
}

@article{mcgrath_silico_2021,
    title = {In silico evaluation and optimisation of magnetic resonance elastography of the liver},
    volume = {66},
    issn = {0031-9155},
    url = {https://dx.doi.org/10.1088/1361-6560/ac3263},
    doi = {10.1088/1361-6560/ac3263},
    language = {en},
    number = {22},
    urldate = {2025-02-25},
    journal = {Physics in Medicine \& Biology},
    author = {McGrath, Deirdre M and Bradley, Christopher R and Francis, Susan T},
    month = nov,
    year = {2021},
    pages = {225005},
}

@article{mcgrath_magnetic_2016,
    title = {Magnetic resonance elastography of the brain: {An} in silico study to determine the influence of cranial anatomy},
    volume = {76},
    issn = {1522-2594},
    shorttitle = {Magnetic resonance elastography of the brain},
    url = {https://onlinelibrary.wiley.com/doi/abs/10.1002/mrm.25881},
    doi = {10.1002/mrm.25881},
    language = {en},
    number = {2},
    urldate = {2025-02-26},
    journal = {Magnetic Resonance in Medicine},
    author = {McGrath, Deirdre M. and Ravikumar, Nishant and Wilkinson, Iain D. and Frangi, Alejandro F. and Taylor, Zeike A.},
    year = {2016},
    pages = {645--662},
}

@article{mcgarry_mapping_2022,
    title = {Mapping heterogenous anisotropic tissue mechanical properties with transverse isotropic nonlinear inversion {MR} elastography},
    volume = {78},
    issn = {1361-8415},
    url = {https://www.sciencedirect.com/science/article/pii/S1361841522000834},
    doi = {10.1016/j.media.2022.102432},
    urldate = {2025-02-26},
    journal = {Medical Image Analysis},
    author = {McGarry, Matthew and Van Houten, Elijah and Sowinski, Damian and Jyoti, Dhrubo and Smith, Daniel R. and Caban-Rivera, Diego A. and McIlvain, Grace and Bayly, Philip and Johnson, Curtis L. and Weaver, John and Paulsen, Keith},
    month = may,
    year = {2022},
    pages = {102432},
}

@article{kwon_shear_2009,
    title = {Shear {Modulus} {Decomposition} {Algorithm} in {Magnetic} {Resonance} {Elastography}},
    volume = {28},
    issn = {1558-254X},
    url = {https://ieeexplore.ieee.org/document/5265288/?arnumber=5265288},
    doi = {10.1109/TMI.2009.2019823},
    number = {10},
    urldate = {2025-02-26},
    journal = {IEEE Transactions on Medical Imaging},
    author = {Kwon, Oh In and Park, Chunjae and Nam, Hyun Soo and Woo, Eung Je and Seo, Jin Keun and Glaser, K. J. and Manduca, A. and Ehman, R. L.},
    month = oct,
    year = {2009},
    pages = {1526--1533},
}

@article{mcgarry_suitability_2015,
    title = {Suitability of poroelastic and viscoelastic mechanical models for high and low frequency {MR} elastography},
    volume = {42},
    issn = {2473-4209},
    url = {https://onlinelibrary.wiley.com/doi/abs/10.1118/1.4905048},
    doi = {10.1118/1.4905048},
    language = {en},
    number = {2},
    urldate = {2025-02-26},
    journal = {Medical Physics},
    author = {McGarry, M. D. J. and Johnson, C. L. and Sutton, B. P. and Georgiadis, J. G. and Van Houten, E. E. W. and Pattison, A. J. and Weaver, J. B. and Paulsen, K. D.},
    year = {2015},
    pages = {947--957},
}

@article{pepin_magnetic_2015,
    title = {Magnetic resonance elastography ({MRE}) in cancer: {Technique}, analysis, and applications},
    volume = {0},
    issn = {0079-6565},
    shorttitle = {Magnetic resonance elastography ({MRE}) in cancer},
    url = {https://www.ncbi.nlm.nih.gov/pmc/articles/PMC4660259/},
    doi = {10.1016/j.pnmrs.2015.06.001},
    urldate = {2025-04-08},
    journal = {Progress in nuclear magnetic resonance spectroscopy},
    author = {Pepin, Kay M. and Ehman, Richard L. and McGee, Kiaran P.},
    month = nov,
    year = {2015},
    pmid = {26592944},
    pmcid = {PMC4660259},
    pages = {32--48},
}

@article{sack_structure-sensitive_2013,
    title = {Structure-sensitive elastography: on the viscoelastic powerlaw behavior of in vivo human tissue in health and disease},
    volume = {9},
    issn = {1744-6848},
    shorttitle = {Structure-sensitive elastography},
    url = {https://pubs.rsc.org/en/content/articlelanding/2013/sm/c3sm50552a},
    doi = {10.1039/C3SM50552A},
    language = {en},
    number = {24},
    urldate = {2025-04-08},
    journal = {Soft Matter},
    author = {Sack, Ingolf and Jöhrens, Korinna and Würfel, Jens and Braun, Jürgen},
    month = may,
    year = {2013},
    pages = {5672--5680},
}

@book{hirsch_magnetic_2017,
    address = {Weinheim},
    title = {Magnetic {Resonance} {Elastography}: {Physical} {Background} {And} {Medical} {Applications}},
    isbn = {978-3-527-34008-8},
    shorttitle = {Magnetic {Resonance} {Elastography}},
    language = {eng},
    publisher = {Wiley-VCH Verlag},
    author = {Hirsch, Sebastian and Braun, Jürgen and Sack, Ingolf},
    year = {2017},
}

@article{meyer_comparison_2022,
    title = {Comparison of inversion methods in MR elastography: {An} open‐access pipeline for processing multifrequency shear‐wave data and demonstration in a phantom, human kidneys, and brain},
    volume = {88},
    issn = {0740-3194, 1522-2594},
    shorttitle = {Comparison of inversion methods in {\textless}span style="font-variant},
    url = {https://onlinelibrary.wiley.com/doi/10.1002/mrm.29320},
    doi = {10.1002/mrm.29320},
    language = {en},
    number = {4},
    urldate = {2025-03-06},
    journal = {Magnetic Resonance in Medicine},
    author = {Meyer, Tom and Marticorena Garcia, Stephan and Tzschätzsch, Heiko and Herthum, Helge and Shahryari, Mehrgan and Stencel, Lisa and Braun, Jürgen and Kalra, Prateek and Kolipaka, Arunark and Sack, Ingolf},
    month = oct,
    year = {2022},
    pages = {1840--1850},
}

@article{arndt_dealii_2023,
    title = {The deal.{II} {Library}, {Version} 9.5},
    volume = {31},
    copyright = {De Gruyter expressly reserves the right to use all content for commercial text and data mining within the meaning of Section 44b of the German Copyright Act.},
    issn = {1569-3953},
    url = {https://www.degruyterbrill.com/document/doi/10.1515/jnma-2023-0089/html},
    doi = {10.1515/jnma-2023-0089},
    language = {en},
    number = {3},
    urldate = {2025-04-08},
    journal = {Journal of Numerical Mathematics},
    author = {Arndt, Daniel and Bangerth, Wolfgang and Bergbauer, Maximilian and Feder, Marco and Fehling, Marc and Heinz, Johannes and Heister, Timo and Heltai, Luca and Kronbichler, Martin and Maier, Matthias and Munch, Peter and Pelteret, Jean-Paul and Turcksin, Bruno and Wells, David and Zampini, Stefano},
    month = sep,
    year = {2023},
    pages = {231--246},
}

@article{Nwmark_method_1959,
    title = {A {Method} of {Computation} for {Structural} {Dynamics}},
    volume = {85},
    issn = {0044-7951, 2690-2427},
    url = {https://ascelibrary.org/doi/10.1061/JMCEA3.0000098},
    doi = {10.1061/JMCEA3.0000098},
    language = {en},
    number = {3},
    urldate = {2025-04-08},
    journal = {Journal of the Engineering Mechanics Division},
    author = {Newmark, Nathan M.},
    month = jul,
    year = {1959},
    pages = {67--94},
}

@article{papazoglou_algebraic_2008,
    title = {Algebraic {Helmholtz} inversion in planar magnetic resonance elastography},
    volume = {53},
    issn = {0031-9155},
    doi = {10.1088/0031-9155/53/12/005},
    language = {eng},
    number = {12},
    journal = {Physics in Medicine and Biology},
    author = {Papazoglou, S. and Hamhaber, U. and Braun, J. and Sack, I.},
    month = jun,
    year = {2008},
    pmid = {18495979},
    pages = {3147--3158},
}

@inproceedings{fike_optimization_2011,
    address = {Honolulu, Hawaii},
    title = {Optimization with {Gradient} and {Hessian} {Information} {Calculated} {Using} {Hyper}-{Dual} {Numbers}},
    isbn = {9781624101458},
    url = {https://arc.aiaa.org/doi/10.2514/6.2011-3807},
    doi = {10.2514/6.2011-3807},
    language = {en},
    urldate = {2025-04-09},
    booktitle = {29th {AIAA} {Applied} {Aerodynamics} {Conference}},
    publisher = {American Institute of Aeronautics and Astronautics},
    author = {Fike, Jeffrey and Jongsma, Sietse and Alonso, Juan and Van Der Weide, Edwin},
    month = jun,
    year = {2011},
}

@article{lai_extensions_2008,
    title = {Extensions of the first and second complex-step derivative approximations},
    volume = {219},
    issn = {0377-0427},
    url = {https://www.sciencedirect.com/science/article/pii/S0377042707004086},
    doi = {10.1016/j.cam.2007.07.026},
    number = {1},
    urldate = {2025-04-09},
    journal = {Journal of Computational and Applied Mathematics},
    author = {Lai, K. -L. and Crassidis, J. L.},
    month = sep,
    year = {2008},
    pages = {276--293},
}

@misc{belponer_reduced_2023,
    title = {Reduced {Lagrange} multiplier approach for non-matching coupled problems in multiscale elasticity},
    url = {http://arxiv.org/abs/2309.06797},
    doi = {10.48550/arXiv.2309.06797},
    urldate = {2025-06-17},
    publisher = {arXiv},
    author = {Belponer, Camilla and Caiazzo, Alfonso and Heltai, Luca},
    month = sep,
    year = {2023},
    note = {arXiv:2309.06797},
}

@article{heltai_zunino_2023,
	author = {Luca Heltai and Paolo Zunino},
	journal = {Mathematical Models and Methods in Applied Sciences},
	number = {12},
	pages = {2425-2462},
	title = {Reduced Lagrange multiplier approach for non-matching coupling of mixed-dimensional domains},
	volume = {33},
	year = {2023}}

@article{heltai_multiscale_2019,
    title = {Multiscale modeling of vascularized tissues via nonmatching immersed methods},
    volume = {35},
    issn = {2040-7947},
    doi = {10.1002/cnm.3264},
    language = {eng},
    number = {12},
    journal = {International Journal for Numerical Methods in Biomedical Engineering},
    author = {Heltai, Luca and Caiazzo, Alfonso},
    month = dec,
    year = {2019},
    pmid = {31508902},
    pages = {e3264},
}

@article{barrett_finite_1986,
    title = {Finite element approximation of the {Dirichlet} problem using the boundary penalty method},
    volume = {49},
    issn = {0945-3245},
    url = {https://doi.org/10.1007/BF01389536},
    doi = {10.1007/BF01389536},
    language = {en},
    number = {4},
    urldate = {2025-07-16},
    journal = {Numerische Mathematik},
    author = {Barrett, John W. and Elliott, Charles M.},
    month = jul,
    year = {1986},
    pages = {343--366},
}

@article{berenger_perfectly_1994,
    title = {A perfectly matched layer for the absorption of electromagnetic waves},
    volume = {114},
    issn = {0021-9991},
    url = {https://www.sciencedirect.com/science/article/pii/S0021999184711594},
    doi = {10.1006/jcph.1994.1159},
    number = {2},
    urldate = {2025-07-16},
    journal = {Journal of Computational Physics},
    author = {Berenger, Jean-Pierre},
    month = oct,
    year = {1994},
    pages = {185--200},
}

@misc{johnson_notes_2021,
    title = {Notes on {Perfectly} {Matched} {Layers} ({PMLs})},
    url = {http://arxiv.org/abs/2108.05348},
    doi = {10.48550/arXiv.2108.05348},
    urldate = {2025-07-16},
    publisher = {arXiv},
    author = {Johnson, Steven G.},
    month = aug,
    year = {2021},
    note = {arXiv:2108.05348},
}

@article{benzaken_constructing_2024,
    title = {Constructing {Nitsche}’s {Method} for {Variational} {Problems}},
    volume = {31},
    issn = {1886-1784},
    url = {https://doi.org/10.1007/s11831-023-09953-6},
    doi = {10.1007/s11831-023-09953-6},
    language = {en},
    number = {4},
    urldate = {2026-01-15},
    journal = {Archives of Computational Methods in Engineering},
    author = {Benzaken, Joseph and Evans, John A. and Tamstorf, Rasmus},
    month = may,
    year = {2024},
    pages = {1867--1896},
}

@article{javili_computational_2020,
    title = {The computational framework for continuum-kinematics-inspired peridynamics},
    volume = {66},
    issn = {1432-0924},
    url = {https://doi.org/10.1007/s00466-020-01885-3},
    doi = {10.1007/s00466-020-01885-3},
    language = {en},
    number = {4},
    urldate = {2026-01-15},
    journal = {Computational Mechanics},
    author = {Javili, A. and Firooz, S. and McBride, A. T. and Steinmann, P.},
    month = oct,
    year = {2020},
    pages = {795--824},
}

@book{hughes_finite_1987,
    address = {Englewood Cliffs, N.J.},
    edition = {1. Dr.},
    title = {The finite element method: linear static and dynamic finite element analysis},
    isbn = {9780133170252 9780133170177},
    shorttitle = {The finite element method},
    language = {eng},
    publisher = {Prentice Hall},
    author = {Hughes, Thomas J. R.},
    year = {1987},
}

@article{mura_analytical_2020,
    title = {An analytical solution to the dispersion-by-inversion problem in magnetic resonance elastography},
    volume = {84},
    issn = {1522-2594},
    doi = {10.1002/mrm.28247},
    language = {eng},
    number = {1},
    journal = {Magnetic Resonance in Medicine},
    author = {Mura, Joaquin and Schrank, Felix and Sack, Ingolf},
    month = jul,
    year = {2020},
    pmid = {32141650},
    pages = {61--71},
}

@article{sack_magnetic_2023,
    title = {Magnetic resonance elastography from fundamental soft-tissue mechanics to diagnostic imaging},
    volume = {5},
    copyright = {2022 Springer Nature Limited},
    issn = {2522-5820},
    url = {https://www.nature.com/articles/s42254-022-00543-2},
    doi = {10.1038/s42254-022-00543-2},
    language = {en},
    number = {1},
    urldate = {2026-01-29},
    journal = {Nature Reviews Physics},
    author = {Sack, Ingolf},
    month = jan,
    year = {2023},
    pages = {25--42},
}
\bibliographystyle{plainnat}

\end{document}